%

\input amstex
\documentstyle{amsppt}
\font\alg=cmss10
\magnification=\magstep1
\def\sipkatau{\longrightarrow_{\kern-10pt\tau}\kern10pt }

\topmatter

\title
The sequential topology on complete Boolean algebras
\endtitle

\author
Bohuslav Balcar, Wieslaw Glowczynski and Thomas Jech
\endauthor

\affil 
The Academy of Sciences of Czech Republic\\
Gdansk University\\
The Pennsylvania State University
\endaffil

\thanks
Supported in part by a grant no. 11904 from AV\v CR (Balcar), and by
the National Science Foundation grant DMS--9401275 and by the National
Research Council COBASE grant (Jech). Glowczynski and Jech are both
grateful for the hospitality
of the Center for Theoretical Study in Prague\endthanks

\address
Mathematical Institute of the Academy of Sciences of Czech Republic,
\v Zitn\'a 25, 115 67 Praha 1, Czech republic (Balcar)
\endaddress
\email balcar\@mbox.cesnet.cz \endemail
\address
Institute of Mathematics, Gdansk University, Wita Stwosza 57,
Gdansk 80-952, Poland (Glowczynski)
\endaddress
\email matwg\@halina.univ.gda.pl \endemail
\address
Department of Mathematics, The Pennsylvania State University,
218 McAllister Bldg., University Park, PA 16802, U.S.A. (Jech)
\endaddress
\email jech\@math.psu.edu \endemail
\abstract
{We investigate the sequential topology $\tau_s$ on a complete
Boolean algebra $B$ determined by algebraically convergent
sequences in $B$. We show the role of weak
distributivity of $B$ in separation axioms for the sequential
topology. The main result is that a necessary and sufficient 
condition for $B$ to carry a strictly positive Maharam submeasure 
is that $B$ is ccc and that the space $(B,\tau_s)$ is Hausdorff.
We also characterize sequential cardinals.}
\endabstract
\endtopmatter

\document

\baselineskip 20pt

\subheading{1. Introduction}
We deal with sequential topologies on complete Boolean algebras from
the point of view of separation axioms. Our motivation comes from the
still open Control Measure Problem of D. Maharam (1947) [Ma]. Maharam
asked wheather every $\sigma$-complete Boolean algebra 
that carries a strictly positive continuous submeasure admits a
$\sigma$-additive measure.

Let us review basic notions and facts concerning Maharam's problem.
More details and further information can be found in Fremlin's paper 
[Fr1].

Let $B$ be a Boolean algebra. A {\it submeasure} on $B$ is
a function $\mu :B\to \bold R^+$ with the properties
\roster
\item"(i)" $\mu(\boldkey 0)=0$,
\item"(ii)" $\mu(a)\leq\mu(b)$ whenever $a\leq b$ (monotone),
\item"(iii)" $\mu(a\vee b)\leq \mu(a)+ \mu(b)$ (subadditive).
\endroster
A submeasure $\mu$ on $B$ is 
\roster
\item"(iv)" {\it exhaustive} if $\lim\mu (a_n)=0$ for every sequence
$\{a_n: n \in \omega \}$ of disjoint elements,
\item"(v)" {\it strictly positive} if $\mu(a)=0$ only if $a=\boldkey0$,
\item"(vi)" a (finitely additive) {\it measure} if for any disjoint 
$a$ 
and $b$, $\mu(a\vee b) = \mu(a)+ \mu(b)$.
\endroster

If $B$ is a $\sigma$-complete algebra a submeasure $\mu$ on $B$ is
called a {\it Maharam submeasure} if it is {\it continuous}, 
i.e. $\lim \mu(a_n)=0$ for every decreasing sequence
$\{a_n: n\in\omega\}$ such that  
$\bigwedge\{a_n: n\in\omega\}=\boldkey0$. 
It is easy to see that a measure on 
a $\sigma$-complete algebra is continuous if and only if it is 
$\sigma$-additive.

We consider the following four classes of Boolean
algebras.

{\alg MBA}: the class of all Boolean algebras that carry a strictly 
positive
finitely additive measure.

{\alg McBA}: the class of all {\it measure} algebras, i.e. 
complete Boolean algebras
that carry a strictly positive $\sigma$-additive measure.

{\alg EBA}: the class of all Boolean algebras that carry 
a strictly positive exhaustive submeasure.

{\alg CcBA}: the class of all complete algebras that carry a strictly
positive continuous submeasure.

\medskip

The diagram below shows the obvious relations between these classes:
\bigskip
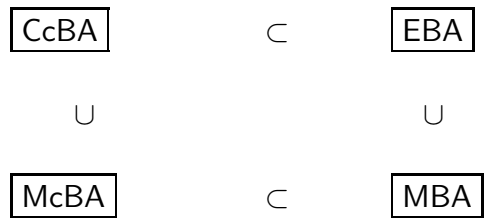
\settabs9\columns
\+&&&{$\boxed{\text{\alg CcBA}}$} &&$\subset$
&{$\boxed{\text{\alg EBA}}$}\cr
\medskip
\+&&&$\qquad\cup$ &&& $\quad\cup$ \cr

\medskip
\+&&&{$\boxed{\text{\alg McBA}}$} &&$\subset$ &{$\boxed{\text{\alg
MBA}}$}\cr

\bigskip

The following theorem whose proof is scattered throughout Fremlin's 
paper
\cite{Fr1} gives additional information. Note that the relations 
between
the classes with measure are the same as between the classes with 
submeasure.

\proclaim{Theorem 1.1} (i) The class {\alg MBA} consists exactly of
all subalgebras of algebras in {\alg McBA}.

(ii) The class {\alg EBA} consists exactly of all subalgebras of 
algebras in
{\alg CcBA}.

(iii) The class {\alg McBA} consists of all algebras in {\alg
MBA} that are complete and weakly distributive.

(iv) The class {\alg CcBA} consists of all algebras in {\alg EBA} that 
are
complete and weakly distributive.
\endproclaim

The problem whether {\alg CcBA} = {\alg McBA} is the problem of
Maharam mentioned above. It follows from Theorem 1.1 that the
problem is equivalent to the problem whether {\alg EBA} = {\alg MBA}.

The class {\alg MBA} is closed under regular completions:
Let $B$ be a Boolean algebra and let $\mu$ be
a finitely additive strictly positive measure. It follows from \cite{Ke} that 
$\mu$ can be extended to a strictly positive measure on the 
completion $\bar B$. 

Similarly, the class {\alg EBA} is closed under regular completions
(this was kindly pointed to us by S. Koppelberg): Let $B$ be a Boolean
algebra and let $\mu$ be a strictly positive exhaustive submeasure.
By \cite{Fr1}, $B$ can be embedded into a complete Boolean algebra $A$
such that $\mu$ can be extended to a strictly positive exhaustive
submeasure on $A.$ By Sikorski's Extension Theorem (\cite{Ko}, p. 70),
the completion $\bar B$ embeds in $A,$ and so $\bar B$ also carries
a strictly positive exhaustive submeasure.

Let us consider an algebra $B\in${\alg CcBA} and let
$\mu$ be a strictly positive Maharam submeasure on $B$. The submeasure
$\mu$ determines a topology on $B$: $(B,\varrho_\mu)$ is a metric 
space with the distance defined by 
$\varrho_\mu(a,b)=\mu(a\vartriangle b)$ for any $a,b\in B$.
If $\nu$ is another such submeasure then  $\varrho_\mu$ and 
$\varrho_\nu$ are
equivalent; they determine the same topology on $B$. In \cite{Ma},
Maharam studied a sequential topology on complete
Boolean algebras from the point of view of metrizability.

We study sequential topologies on complete Boolean algebras in a 
more general setting.
Our goal is to show that the sequential topology $\tau_s$ on a ccc 
complete 
Boolean algebra $B$ is Hausdorff
if and only if $B$ carries a strictly positive Maharam submeasure. 
Following \cite{AnCh} and \cite {Pl} we say that a cardinal $\kappa$ 
is a 
{\it sequential cardinal} if there
exists a continuous real-valued function on the space $(\Cal
P(\kappa),\tau_s)$
which is not continuous with respect to the product topology. 
We prove that $\kappa$ is a sequential cardinal if and only if 
$\kappa$ is uncountable
and there is a nontrivial Maharam submeasure on the algebra $\Cal
P(\kappa)$.

\subheading{2. Sequential topology}

Let us review some notions from topology.

{\bf 2.1. Definition.} Let $(X,\tau)$ be a topological space. The 
space
$X$ is
\roster
\item"(i)" {\it sequential} if a subset $A\subseteq X$ is closed
whenever it contains all limits of $\tau$-convergent sequences of
elements of $A$;
\item"(ii)" {\it Fr\'echet} if for every $A\subseteq X$
$$\text{cl}_\tau(A)=\{x\in X:(\exists \langle x_n:n\in
\omega\rangle\subseteq A)\ x_n\sipkatau x\}.$$
\endroster

It is clear that every Fr\'echet space is sequential.

Now, consider a complete Boolean algebra $B$; $\sigma$-completness is
sufficient for the following definition. For a sequence
$\langle b_n: n\in \omega\rangle$ of elements of $B$ we denote

$$\overline{\lim}\, b_n=
\bigwedge\limits_{k\in\omega}\bigvee\limits_{n\ge
k}b_n\qquad\text{and} \qquad\underline{\lim}\,
b_n=\bigvee\limits_{k\in\omega}\bigwedge\limits_{n\ge k}b_n\ .$$

We say that a sequence  $\langle b_n\rangle$ {\it algebraically 
converges}
to an element
$b\in B$, $b_n\longrightarrow b$,
 if $\overline{\lim}\, b_n=\underline{\lim}\, b_n=b$.

A sequence $\langle b_n\rangle$ algebraically converges
if and only if there exist an increasing sequence $\langle a_n\rangle$
and a decreasing sequence $\langle c_n\rangle$ such that
$a_n\leq b_n\leq c_n$ for all $n\in \omega,$ and
$\bigvee\limits_{n\in\omega}a_n=\bigwedge\limits_{n\in\omega}c_n$.

{\bf 2.2.} We summarize basic properties of $\longrightarrow $:
\roster
\item "(i)" Every sequence has at most one limit;
\item "(ii)" for a constant sequence $\langle x:n\in\omega\rangle$, we
have  $\langle x:n\in\omega\rangle
\longrightarrow x$;
\item "(iii)" $x_n\longrightarrow\boldkey 0$ iff $\overline{\lim}\,
x_n=\boldkey 0$;
\item "(iv)" if the $x_n$'s are pairwise disjoint then 
$x_n\longrightarrow \boldkey 0$;
\item "(v)" $\overline{\lim}\, (x_n \vee y_n) = \overline{\lim}\, x_n
\vee \overline{\lim}\, y_n$;
\item "(vi)" if $x_n\longrightarrow x$ and $y_n\longrightarrow y$ then
$x_n \vee y_n\longrightarrow x\vee y$ and $-x_n\longrightarrow -x$;
\item "(vi)" if $\langle x_n\rangle$ is increasing then
$x_n\longrightarrow \bigvee\limits_{n\in\omega}x_n$.
\endroster

\subheading{2.3. Sequential topology on $B$} Let us consider all 
topologies
$\tau$ on $B$ with the following property:
$$\text{ if } x_n\longrightarrow x\  \text{ then }  x_n\sipkatau x.$$
There is a largest topology with respect to the inclusion among 
all such topologies.
We denote it  $\tau_s$ and call it the {\it sequential topology} on
$B$.

The topology $\tau_s$ can be described as follows, by definining the
closure operation: For any subset $A$ of the algebra $B$ let 
$$
u(A) = \{x : x \text{ is the limit of a sequence } \{x_n\} \text
{ of elements of } A\}.
$$
The closure of a set $A$ in the topology $\tau_s$ is obtained by an 
iteration
of $u$:
$$
cl_{\tau_s} (A) = \bigcup_{\alpha <\omega_1} u^{(\alpha)}(A),
$$
where $u^{(\alpha+1)}(A) = u(u^{(\alpha)}(A))$, and for a limit 
$\alpha$, $u^{(\alpha)}(A) = \bigcup_{\beta<\alpha} u^{(\beta)}$.

It is clear that the topology $\tau_s$ is $T_1$, i. e. every singleton 
is a 
closed set. Moreover, $(B, \tau_s)$ is a Fr\'echet space if and only 
if
$cl(A)=u(A)$ for every $A\subseteq B.$

We remark that a sequence $\{x_n\}$ converges to $x$ topologically if 
and
only if every subsequence of $\{x_n\}$ has a subsequence that 
converges
to $x$ algebraically.

\proclaim{Example: Measure algebras}\endproclaim
Let $B$ be a complete Boolean algebra carrying a strictly positive 
$\sigma-$additive
measure $\mu.$ 
For any $a,b\in B,$ let
$$
\rho(a,b) =\mu(a\vartriangle b);
$$
$\rho$ is a metric on $B$ and the topology given by $\rho$ coincides 
with the
sequential topology. Hence $(B,\tau_s)$ is metrizable.

\medpagebreak

Maharam's Control Measure Problem is equivalent to the question
whether there exist complete Boolean algebras other
than the algebras in the class {\alg McBA} 
for which the sequential topology is metrizable.

\subheading {Properties of the topology $\tau_s$}

\proclaim
{Proposition 2.4} (i) The operation of complement is continuous (and 
hence
a homeomorphism).

(ii) For a fixed $a$, the function $a\lor x$ is a continuous function 
of $x.$

(iii) For a fixed $a$, the function $a\vartriangle x$ is a continuous 
function.
\endproclaim

The operation $\lor$ is generally not a continuous function of two 
variables.
As a consequence of (iii), the space $(B,\tau_s)$ is homogeneous: 
given
$a, b\in B,$ there exists a homeomorphism $f$ such that $f(a)=b$, 
namely
$f(x)=(x\vartriangle b)\vartriangle a).$ The topology $\tau_s$ is 
determined 
by the
family $\Cal N_0$ of all neighborhoods of $0$ as for every $a\in B$ 
and every set $W$, $W$ is a neighborhood of $a$ if and only if 
$a\vartriangle W\in \Cal N_0$.

As a consequence of homogeneity of $(B,\tau_s)$, $B$ does not have 
isolated
points unless $B$ is finite.

\proclaim
{Lemma 2.5} Let $B$ be a $\sigma-$complete algebra. Let $\{u_n\}_{n=
0}^\infty$
be an antichain in $B$, and let $U$ be a neighborhood of $\,\boldkey0$. 
Then there exists
a $k$ such that $B\restriction\bigvee_{n\ge k} u_n \subset U.$
\endproclaim

\demo
{Proof} If not, then for every $k$ there exists an element $x_k$ below
$\bigvee_{n\ge k} u_n$ such that $x_k\notin U.$ But then the sequence
$\{x_k\}_k$ converges to $0$ and so, because $U\in \Cal N_0,$ there 
exists 
some $k_0$ such that $x_k\in U$ for all $k>k_0.$ A contradiction.
\qed
\enddemo

A subset $D$ of a Boolean algebra $B$ is {\it dense} if for every $b
\in B,
b\ne 0$ there is some $d\in D, d\ne 0$ such that $d\le b$. $D$ is {\it 
downward
closed} if $a<d\in D$ implies $a\in D$. 

If $H$ is a downward closed subset of $B$ then $H\vartriangle H=
H\lor H,$ and so if $H$ is also an open set then so is $H\lor H.$ 

A downward closed dense set is called
{\it open dense}. Since we consider a topology on $B$ we shall call 
dense
and open dense set {\it algebraically dense} and {\it algebraically 
open
dense} to avoid confusion with the corresponding topological terms. 

\proclaim
{Corollary 2.6} (i) Every neighborhood of $\boldkey 0$ contains all 
but finitely many 
atoms.

(ii) If $B$ is atomless then every neighborhood of $\boldkey 0$ 
contains
an algebraically open dense subset of $B$.

(iii) If $B$ is atomless and ccc, then for every $U\in \Cal N_0$ there 
exists
a $k$ such that $\boldkey 1\in U \vartriangle \dots \vartriangle U$ 
($k$ times).
\endproclaim

\demo
{Proof of (ii)} Let $V$ be a neighborhood of $\boldkey 0.$ If $V$ does 
not contain
an algebraically open dense set then $B-V$ is algebraically dense 
below 
some $u\ne 0$ and hence contains
a pairwise disjoint set $\{x_n\}_n.$ But then $\lim x_n = 0$ and so 
there is 
some $n$ such that $x_n\in V;$ a contradiction.

(iii) As $U$ is algebraically dense in $B$, 
there exists a maximal antichain
of $B$ included in $U$, and by ccc the antichain is countable:
$\{u_n\}_n \subset U.$ There exists a $k$ so that $u=\bigvee_{n> k} u_n
\in U,$
and then $u_0 \lor u_1 \lor \dots \lor u_k \lor u =1.$
\qed\enddemo

\proclaim
{Proposition 2.7}
If $B$ is atomless and ccc, then $(B,\tau_s)$ is connected.
\endproclaim

\demo
{Proof} Assume that there are two disjoint nonempty clopen sets $X$ 
and $Y$ 
with $X\cup Y =B$
such that $\boldkey 0\in X$, and let $a\in Y$ be arbitrary. Let $C$ be 
a maximal
chain in $B$ such that $\inf C = 0$ and $\sup C = a.$ Let $x = \sup 
(C\cap X);$
by ccc, $x$ is the limit of a sequence in $C\cap X$ and therefore $x
\in X.$
Let $y = \inf (Y\cap \{c\in C: c\ge x\}).$ Using the ccc again we have 
$y\in Y ,$ and clearly $x< y.$ By maximality, both $x$ and $y$ are in 
$C.$
Since $B$ is atomless, there exists some
$z$ with $x<z<y.$ This contradicts the maximality of $C.$
\qed\enddemo

\proclaim
{Lemma 2.8} (i) An ideal $I$ on a $\sigma-$complete Boolean algebra
$B$ is a closed set in the sequential topology
if and only if it is a $\sigma-$complete ideal.

(ii) If $I$ is a $\sigma-$ideal on $B$ then the sequential topology on
the quotient algebra $B/I$ is the quotient topology of $\tau_s$ given 
by the
canonical projection.

(iii) In (ii), the canonical projection onto $B/I$ is a closed 
mapping.

(iv) If $\tau_s$ is Fr\'echet then so is the quotient topology.
\endproclaim

\subheading{3. Fr\'echet spaces}

We shall now consider those complete Boolean algebras for which the 
sequential
topology is Fr\'echet. We will show that this is equivalent to an 
algebraic
property. First we make the following observation:

\proclaim{Proposition 3.1} If $(B,\tau_s)$ is a Fr\'echet space then 
for
every $V\in \Cal N_0$ there is some $U\subseteq V$ in $\Cal N_0$ such 
that
$U$ is downward closed.
\endproclaim

\demo{Proof} If $V\in \Cal N_0$, consider the set
$$
X= \{a\in B: \text{ there exists some } b\le a \text{ such that } b
\notin V\},
$$ 
and let $u(X)$ be the set of all limits of sequences in $X.$ As 
$\tau_s$ is Fr\'echet,
$u(X)$ is the closure of $X.$ We shall prove that the set $U=B-u(X)$ 
is downward closed and contains $\boldkey 0.$

For the first claim it suffices to show that $a\in u(X)$ and $a<b$ 
implies 
$b\in u(X).$ Thus let $a=\lim\,a_n$ with $a_n\in X.$ It follows that
$b=\lim\,(a_n \lor b),$ and since $a_n\lor b \in X, $ we have $b\in 
u(X).$

To see that $\boldkey 0\notin U,$ assume that $\{a_n\} \subseteq X$ 
and $\lim\,a_n =0.$
Then there are
$x_n\le a_n$ in $B-V$, but this is impossible because $\lim\,x_n=0.$ 
Hence $\boldkey 0$ is not in $u(X).$
\qed
\enddemo

Thus if $(B,\tau_s)$ is Fr\'echet, its topology is determined by the 
set
$\Cal N^d_0$ of all
$U\in \Cal N_0$ that are downward closed. $\Cal N_0^d$ is a 
neighborhood base
of $\boldkey 0.$

\proclaim
{Definition 3.2} Let $\kappa$ be an infinite cardinal. A Boolean 
algebra $B$ is 
$(\omega,\kappa)$-weakly distributive if for every sequence $\{P_n\}$ 
of 
maximal antichains, each of size at most $\kappa,$ there exists a 
dense set
$Q$ with the property that each $q\in Q$ meets only finitely many 
elements
of each $P_n.$ $B$ is {\it weakly distributive} if it is $(\omega,
\omega)$-weakly distributive.
\endproclaim

If $B$ is a $\kappa^+$-complete Boolean algebra then $B$ is $(\omega,
\kappa)$-weakly distributive if and only if it satisfies the following
distributive law:
$$
\bigwedge_n \bigvee_\alpha a_{n \alpha} = 
\bigvee_{f:\omega\to [\kappa]^{<\omega}} \bigwedge_n \bigvee_{\alpha\in 
f(n)} a_{n \alpha}.
$$

We recall two frequently used cardinal characteristics.
\proclaim {Definition 3.3} 
The {\it splitting number} is the least cardinal $\bold s$ of a family 
$\Cal S$
of infinite subsets of $\omega$ such that for every infinite 
$X\subseteq\omega$
there is some $S\in\Cal S$ such that both $X\cap S$ and $X-S$ are 
infinite.
($S$ ``splits'' $X$.)

The {\it bounding number} is the least cardinal $\bold b$ of a family 
$\Cal F$
of functions from $\omega$ to $\omega$ such that $\Cal F$ is 
unbounded; i.e.
for every $g\in \omega^\omega$ there is some $f\in \Cal F$ such that 
$g(n)
\le f(n)$ for infinitely many $n.$
\endproclaim

\medpagebreak

The following characterization of Fr\'echet spaces $(B,\tau_s)$ uses 
the
cardinal invariant $\bold b$ and is similar to several other results
using $\bold b$ such as \cite{BlJe}.
A consequence of Theorem 3.4 is that
$(P(\kappa),\tau_s)$ is a Fr\'echet space if and only if $\kappa < 
\bold b$.

\proclaim{Theorem 3.4} Let $B$ be a complete Boolean algebra.
The sequential space $(B,\tau_s)$ is Fr\'echet if and
only if $B$ is weakly distributive and satisfies
the $\bold b$-chain condition.
\endproclaim

We shall first reformulate the condition stated in Theorem 3.4. 
Let $B$ be a complete Boolean algebra. We shall call a matrix
$\{a_{mn}\}$ {\it increasing} if each row $\{a_{mn}:n\in \omega\}$ is 
an
increasing sequence with limit $\boldkey 1.$ Note that $B$ is weakly 
distributive if
and only if for every
increasing matrix $\{a_{mn}\}$,
$$
\bigvee_{f\in\omega^\omega} \underline{\lim}\,a_{m,f(m)} = 1.
$$

\proclaim{Lemma 3.5} A complete Boolean algebra $B$ is weakly 
distributive
and satisfies the $\bold b$-chain condition if and only if for every 
increasing
matrix $\{a_{mn}\}$ there exists a function $f\in\omega^\omega$ such 
that
$\lim a_{m,f(m)} = 1.$
\endproclaim

\demo{Proof} First let $B$ be weakly distributive and $\bold b$-c.c., 
and let
$\{a_{mn}\}$ be an increasing matrix. By the $\bold b$-chain condition
there exists a set $F\subset 
\omega^\omega$ of size less than $\bold b$ such that $\bigvee_{f\in F}
\underline{\lim}\,a_{m,f(m)} =1.$
Let $g:\omega \to\omega$ be an upper bound of $F$ under eventual 
domination.
Since the matrix is increasing, we have $\underline{\lim}\,a_{m,f(m)}
\le \underline{\lim}\,a_{m,g(m)}$ for every $f\in F$. Therefore $\lim 
a_{m,g(m)}
=1.$

Conversely, assume that the condition holds. Then $B$ is weakly 
distributive,
and we shall verify the $\bold b$-chain condition. Thus let $W$ be a 
partition
of $\boldkey 1$; we shall prove that $|W|<\bold b.$ Let 
$\{f_u:u\in W\}$ be any family
of functions from $\omega$ to $\omega$ indexed by elements of $W.$ For 
each
$m$ and each $n$ we let
$$
a_{mn} = \bigvee \{u\in W: f_u(m)<n\}.
$$
The matrix $\{a_{mn}\}$ is increasing and therefore there exists a 
function
$g:\omega\to\omega$ such that $\lim\,a_{m,g(m)} = 1.$ Since $W$ is an 
antichain,
it follows that for
any $u\in W$ there is some $m_u$ such that $u\le a_{m,g(m)}$ for every
$m\ge m_u.$ Hence $f_u (m) < g(m)$ for every $m\ge m_u$ and it follows 
that
$g$ is an upper bound of the family $\{f_u:u\in W\}.$ Therefore every
family of functions of size $|W|$ is bounded and so $|W|<\bold b.$
\qed
\enddemo

\demo{Proof of Theorem 3.4} We wish to show that the condition in 
Lemma 3.5
is necessary and sufficient for the space $(B,\tau_s)$ to be 
Fr\'echet.
To see that the condition holds if $(B,\tau_s)$ is Fr\'echet, we 
recall
\cite{Ma} that for $(B,\tau_s),$ being Fr\'echet
is equivalent to the following statement:
whenever $\{x_{mn}\}$,  $\{y_m\}$ and $z$ are such that $\lim_n 
x_{mn}=y_m$
for each $m$ and $\lim_m y_m = z$, then there is an $f:\omega\to
\omega$ 
such that $\lim_m x_{m,f(m)} = z.$

To show that the condition implies that $(B,\tau_s)$ is Fr\'echet,
let $\{x_{mn}\}, \{y_m\}$ and $z$
be as above. For each $m$ and each $n$ let $u_{mn} = x_{mn}\vartriangle
(-y_m),$ and let $a_{mn} = \bigwedge_{k\ge n} u_{mk}.$ For each $m$, 
$\lim_n u_{mn}=1;$ the matrix $\{a_{mn}\}$ is increasing, with each 
row 
converging to $\boldkey 1$ and so
there exists some $f:\omega\to\omega$ such that $\lim \,a_{m,f(m)} = 
1.$
It follows that $\lim_m \bigwedge_{k\ge f(m)} u_{mk} =1,$ and so
$\lim\,(x_{m,f(m)}\vartriangle (-z) = 
\lim\,(x_{m,f(m)}\vartriangle(-y_m))=\lim \, u_{m,f(m))}= 1.$ 
Hence $\lim x_{m,f(m)} = z.$
\qed\enddemo

We conclude with the following observation that we shall use in Section 5.

\proclaim{Lemma 3.6} (a) For every set $A\subseteq B$, $cl(A)=
\bigcap\{A\vartriangle V: V\in \Cal N_0\}.$

(b) If $(B,\tau_s)$ is Fr\'echet and $A$ is downward closed then 
$cl(A)=
\bigcap\{A\vee V: V\in\Cal N^d_0\}$, and $cl(A)$ is downward closed.
\endproclaim

\demo{Proof} (a) For any $x\in B,$ $x\in cl(A)$ iff for all $V\in \Cal 
N_0,$
$(V\vartriangle x) \cap A \ne \emptyset$, i.e. there exist $v\in V$ 
and $a\in A$
such that $v\vartriangle x = a.$ The latter is equivalent to 
$x=a\vartriangle v$, or $x\in A\vartriangle V.$ 

(b) If both $A$ and $V$ are downward closed then 
$A\vee V = A\vartriangle V.$
\qed\enddemo

\proclaim{Corollary 3.7} For every $U\in\Cal N_0,$ 
$cl(U)\subseteq U\vartriangle
U.$ If $(B,\tau_s)$ is Fr\'echet, then for every
$U\in\Cal N_0^d,$ $cl(U) \subseteq U\vee U.$
\endproclaim

\subheading {4. Separation axioms}

We will now discuss separation axioms for the topology $\tau_s$. 
We see immediately that the sequential topology on $B$ is $T_1.$ The 
space is 
Hausdorff if and only if every point
$b\ne 0$ can be separated from $\boldkey 0$ which is equivalent to the 
statement
that for every $b\ne 0$ there exists some $V\in \Cal N_0$ such that
$b\notin V\vartriangle V$.

\proclaim
{Theorem 4.1} If $(B,\tau_s)$ is a Hausdorff space then $B$ is 
$(\omega,\omega_1)-$weakly distributive.
\endproclaim

We first prove a weaker statement, namely that Hausdorff implies weak
distributivity:

\proclaim
{Lemma 4.2} If $B$ is not weakly distributive then there exists an $a\ne 0$
such that $c\in cl(U)$ for every $c\le a$ and every $U\in \Cal N_0.$ 
Hence $(B,\tau_s)$ is not Hausdorff.
\endproclaim

\demo {Proof} Assume that $B$ is not
$(\omega,\omega)-$weakly distributive. There is some $a\ne 0$ and
there exists an infinite matrix $\{a_{m n}\}$
such that each row  is a partition of $a$, and for any nonzero $x\le a$
there is some $m$ such that $x\wedge a_{m n} \ne 0$ for infinitely many $n.$ 

Let $c\le a$ and let $U$  be an arbitrary neighborhood of $\boldkey 0.$ We will
show that $c\in cl(U).$ For every $m$ and every $n$ let
$y_{m n} = c \wedge\bigvee_{i\ge n} a_{m i}.$ Since the sequence $\{y_{0 n}\}$
converges to $\boldkey 0$ there exists some $n_0$ such that $y_{0 
n_0}\in U;$
let $x_0 = y_{0 n_0}.$ Next we consider the sequence $\{y_{1 n}\lor x_0\}.$
This sequence converges to $x_0$ and so there exists some $n_1$ such that
$x_1\in U$ where $x_1= y_{1 n_1}\lor x_0.$ We proceed by induction and
obtain a sequence $\{n_m\}$ and 
an increasing sequence $\{x_m\}$ of elements of $U.$ This sequence
converges to $c$ because otherwise if we let $b\ne 0$ be the complement of
$\bigvee_n x_n$ in $c$ then $b\leq \bigwedge_m \bigvee_{i<n_m} a_{mi}$ and
so $b$ meets only finitely many elements in each row of the matrix.
Hence $c\in cl(U).$
\qed\enddemo

\demo{Proof of Theorem 4.1} Let $(B,\tau_s)$ be a Hausdorff space.
To prove that $B$ is $(\omega,\omega_1)-$weakly distributive, let
$$
A=\{a_{n\alpha} : n\in\omega,\alpha\in\omega_1\}
$$
be a matrix such that each row is a partition of $\boldkey 1.$ Denote 
$X$ the set
of all those $x\in B$ that meet at most countably many elements of 
each row 
of $A$. As $B$ is $(\omega,\omega)-$weakly distributive, for every 
nonzero 
$x\in X$ there is a nonzero $y\le x$ that meets only finitely many 
elements 
of each row of $A.$ Thus we complete the proof by showing that 
$\bigvee X=1.$

Assume otherwise; without loss of generality we may assume that every 
$x\ne 0$
meets uncountably many elements of at least one row of $A.$ Then the 
matrix $A$
represents a Boolean-valued name for a cofinal function from $\omega$ 
into
$\omega_1.$ Thus $B$ collapses $\omega_1$ and therefore there exists a 
matrix
$$
\{b_{n\alpha} : n\in \omega, \alpha\in\omega_1\}
$$
such that each row and each column is a partition of $\boldkey 1$ (the 
name for a
one-to-one mapping of $\omega$ onto $\omega_1$). We get a 
contradiction to
Hausdorffness by showing that $\boldkey 1$ is in the closure of every 
$V\in \Cal N_0.$

Let $V\in \Cal N_0$ be arbitrary. By Lemma 2.5 there is for every 
$\alpha\in
\omega_1$ some $n_\alpha \in\omega$ such that $v_\alpha = 
\bigvee_{i\ge n_\alpha} b_{i\alpha} \in V.$ Thus there exists some $n$, 
and an infinite set 
$\{\alpha_k\}_k$ such that $n_{\alpha_k} = n$ for all $k.$ Now, by 2.2 (v),
$$
\overline{\lim_k} \bigvee_{i<n} b_{i\alpha_k} =
\bigvee_{i<n} \overline{\lim_k}\, b_{i\alpha_k} = 0.
$$
Therefore $\lim_k v_{\alpha_k} = 
\boldkey 1$ and so $\boldkey 1$ is in the closure of $V.$
\qed\enddemo

Concerning $(\omega,\omega_1)-$weak distributivity we refer to Namba's
work \cite{Na} which shows that it may or may not be equivalent to
$(\omega,\omega)-$weak distributivity. If $\bold b = \omega_1$ then
$(\omega,\omega)-$weak distributivity and $(\omega,\omega_1)-$weak 
distributivity are equivalent, and there is a model of ZFC in which
they are not equivalent. Below (4.5.i) we give another example of
a complete Boolean algebra 
that is $(\omega,\omega)-$weakly distributive but not $(\omega,
\omega_1)-$weakly distributive.

Theorem 4.1 cannot be extended by replacing $\omega_1$ by $\infty$:
Example 4.5.ii, due to Prikry \cite{Pr}, provides a complete Boolean
algebra that is Hausdorff (therefore weakly distributive) but not
$(\omega,\kappa)-$weakly distributive, for a measurable $\kappa.$

In view of Theorems 3.4 and 4.1 the question arises about the relative
strength of being a Hausdorff space and being a Fr\'echet space.
Example 4.3 below shows that Hausdorff does not imply Fr\'echet:
the space $(P(\bold b),\tau_s)$ is Hausdorff but not Fr\'echet.

For the other direction, see Examples 4.4 and 4.5. If $T$ is a Suslin 
tree
then $(B(T),\tau_s)$ is Fr\'echet but not Hausdorff. 

\proclaim {4.3. Example}\endproclaim 
For every infinite cardinal $\kappa$ the space
$(P(\kappa),\tau_s)$ is Hausdorff. This is because each principal 
ultrafilter
on $\kappa$ is a closed and open subset of $P(\kappa)$. 

We identify $P(\kappa)$ with $2^\kappa$ (via characteristic 
functions).
For each $\alpha\in \kappa$ the sets $\{X\subseteq\kappa : \alpha\in 
X\}$
and its complement $\{X\subseteq\kappa : \alpha\notin X\}$ are closed 
under limits of
sequences and so are both closed and open. This implies that the 
topology
$\tau_s$ extends the product topology, and the space $(P(\kappa),
\tau_s)$
is totally disconnected Hausdorff space. If $\kappa=\aleph_0$ then 
$\tau_s$
is equal to the product topology. To see this, let 
$U\subseteq P(\omega)$ be an
open set in the sequential topology and let $A\in U$. For each $n$ let 
$S_n$
denote the basic open set (in the product topology) 
$\{X\subseteq\omega :
X\cap n = A\cap n\}$. It suffices to show that  $U$ contains some
$S_n$ as a subset. If not, there exists for each $n$ some $X_n \in S_n 
- U.$
But $A = \lim_n X_n$, and since the complement of $U$ is closed, 
$A\notin U$;
a contradiction.

When $\kappa$ is an uncountable cardinal, the space $(P(\kappa),
\tau_s)$ is not
compact and so $\tau_s$ is strictly stronger than the product 
topology.

By \cite{Tr} the space $(P(\kappa),\tau_s)$ is sequentially compact if 
and
only if $\kappa < \bold s$, the splitting number.

By \cite{Gl}, $(P(\kappa),\tau_s)$ is regular if and only if $\kappa=
\omega.$
See Corollary 4.7.

\proclaim{4.4. Example: Aronszajn trees}\endproclaim
Let $T$ be an Aronszajn tree and assume that each node has at least 
two
immediate successors. Let $B(T)$ denote the complete Boolean algebra 
that
has upside down $T$ as a dense set. We will show that $(B(T),\tau_s)$ 
is
not a Hausdorff space. This shows that the converse of Theorem 4.1 is 
not
provable: if $T$ is a Suslin tree then $B(T)$ is a ccc $\omega$-
distributive
Boolean algebra.

We prove that $\boldkey 0$ and $\boldkey 1$ cannot be separated by 
open sets: we show that
for every open neighborhood $V$ of $\boldkey 0$, 
$\boldkey 1\in V\vartriangle V.$ Let 
$V\in\Cal N_0.$ For every $\alpha\in \omega_1,$ the $\alpha$th level 
$T_\alpha$
of the tree is a countable partition of $\boldkey 1$ and so there 
exists a finite
set $u_\alpha\subseteq T_\alpha$ such that 
$x_\alpha = \bigvee(T_\alpha - u_\alpha)\in V.$ Let $y_\alpha = 
\bigvee u_\alpha.$ We claim that there is a $\beta$ such that 
$y_\beta\in V;$
this will complete the proof as $\boldkey 1=x_\beta\vartriangle 
y_\beta \in
V\vartriangle V.$

Let $f : [\omega_1]^2 \to \{0,1\}$ be the function defined as follows:
$f(\alpha,\beta)=0$ if $y_\alpha\wedge y_\beta=0$ and $f(\alpha,
\beta)=1$
otherwise. By the Dushnik-Miller Theorem there exists a set 
$I\subseteq
\omega_1$, either homogeneous in color $0$ and of size $\aleph_0$, or
homogeneous in color $1$ and of size $\aleph_1.$ The latter case is 
impossible
because the $u_\alpha$s are finite sets in an Aronszajn tree (see 
\cite{Je}, Lemma 24.2). Hence there is an infinite set $\{\alpha_n:n
\in\omega\}$
such that the $y_{\alpha_n}$ are pairwise disjoint. Thus the sequence
$\{y_{\alpha_n}\}$ converges to $\boldkey 0$ and so there exists some 
$n$ such 
that $y_{\alpha_n}\in V.$

\proclaim {4.5. Examples using large cardinals}\endproclaim 
(i) Assume that there exists a nontrivial $\aleph_2$-saturated
$\sigma$-ideal $I$ on $P(\omega_1)$, and assume that $\bold b =
\aleph_2.$
Both these assumptions are consequences of Martin's Maximum (MM), with
$I=$ the nonstationary ideal.

Let $B= P(\omega_1)/I$. B is a complete Boolean algebra and satisfies
the $\aleph_2-$chain condition.
Since $\bold b=\aleph_2,$ the space $(P(\omega_1),\tau_s)$ is 
Fr\'echet,
and so by Lemma 2.8, $(B,\tau_s)$ is Fr\'echet. Therefore $B$ is 
weakly
distributive.

Since forcing with $B$ collapses $\aleph_1,$ $B$ is not 
$(\omega,\omega_1)-$weakly distributive, and hence $(B,\tau_s)$ is not
Hausdorff.

The space $(P(\omega_1),\tau_s)$ is separable: this follows from MM,
specifically from $\bold p=\aleph_2$, cf. \cite{Fr0}, \cite{To} and 
\cite{Ro}.
Hence $(B,\tau_s)$ is separable, and so the complete Boolean algebra 
$B$ is
countably generated.

This example is in the spirit of \cite{Gl} where a similar example is
presented using MA and a measurable cardinal.

\medpagebreak

(ii) Let $\kappa$ be a measurable cardinal, and let $B$ be the 
complete
Boolean algebra associated with Prikry forcing. $B$ is not
$(\omega,\kappa)-$weakly distributive as it changes the cofinality of
$\kappa$ to $\omega.$ But the space $(B,\tau_s)$ is Hausdorff: For
any $a\in B^+$ there is a $\kappa-$complete ultrafilter on $B$ 
containing
$a,$ cf. \cite{Pr}. Every such ultrafilter is a clopen set in 
$(B,\tau_s).$

Thus being Hausdorff does not imply $(\omega,\infty)-$weak 
distributivity
of $(B,\tau_s).$

\medpagebreak

A topological space is {\it regular} if points can be separated from 
closed sets;
equivalently, for every point $x$ and its neighborhood $U$ there 
exists an
open set $V$ such that $x\in V$ and $cl(V)\subseteq U.$ The space
$(B,\tau_s)$ is regular if and only if for every $U\in \Cal N_0$ there 
is
some $V\in\Cal N_0$ such that $cl(V)\subseteq U.$

A result proved independently in \cite{Tr} and \cite{Gl} states that
the atomic algebra $P(\omega_1)$ is not regular. The following lemma 
uses the 
method employed in these papers.

\proclaim{Lemma 4.6} In the space $(P(\omega_1),\tau_s)$ for every 
$V\in
\Cal N_0$ there exists a closed unbounded set $C\subset\omega_1$ such 
that
for every $\beta\in C,$ $\omega_1 - \beta \in cl(V).$
\endproclaim

\demo{Proof} Let $V$ be an open neighborhood of $\emptyset.$
Let $\{A_{\alpha n}: \alpha\in\omega_1,n\in\omega\}$ be an Ulam 
matrix,
i.e. a double array of subsets of $\omega_1$ with the following 
properties:
$$
A_{\alpha n}\cap A_{\alpha m}=\emptyset\qquad (n\ne m),
$$
$$
A_{\alpha n}\cap A_{\beta n}=\emptyset\qquad (\alpha\ne\beta),
$$
$$
\bigcup_{n\in\omega} A_{\alpha n} = \omega_1 - \alpha.
$$
By lemma 2.5 there exists for each $\alpha$ some $k_\alpha$ such that
$X_\alpha = \bigcup_{n\ge k_\alpha} A_{\alpha n}$ is in $V.$
There exist some $k$ and an uncountable set $W$ such that 
$k_\alpha = k$ for every $\alpha\in W.$ Let $C$ be the set of
all limits of increasing sequences  of ordinals in $W.$
We claim that for every $\beta\in C,$ $\omega_1 - \beta \in cl(V).$

Let $\alpha_0 < \alpha_1 < \dots < \alpha_n < \dots $ be in $W$ such 
that
$\beta=\lim_n\,\alpha_n.$
Note that $\overline{\lim}_n \bigcup_{i<k} A_{\alpha_n i}$ $=
\bigcup_{i<k} \overline{\lim}_n A_{\alpha_n i} = \emptyset,$
and hence $X = \lim_n X_{\alpha_n} = \omega_1 -\beta.$
Therefore   $X\in cl(V).$ 
\qed\enddemo

\proclaim{Corollary 4.7} The space $(P(\omega_1),\tau_s)$ is not 
regular.
\endproclaim

\demo{Proof} Let $U$ be the set of all $x\subset\omega_1$ whose 
complement is
uncountable. $U$ is an open neighborhood of $\emptyset$ and by Lemma 
4.6,
does not contain $cl(V)$ for any $V\in\Cal N_0.$
\qed\enddemo

\proclaim{Corollary 4.8} If a complete Boolean algebra $B$ does not
satisfy the countable chain condition 
then $(B,\tau_s)$ is not regular.
\endproclaim

\demo{Proof} $B$ contains $P(\omega_1)$ as a complete subalgebra, 
therefore
as a closed subspace. Hence it is not regular.
\qed\enddemo

\proclaim{Corollary 4.9} Let $B=P(\omega_1),$ or more generally, let 
$B$
be a complete Boolean algebra that does not satisfy the countable 
chain
condition. If $\{U_n\}_n$ is a countable subset of $\Cal N_0$ then
$\bigcap_n cl(U_n)$ is uncountable.
\endproclaim

\demo{Proof} This follows easily from Lemma 4.6 when $B=P(\omega_1).$ 
In the general case, $(B,\tau_s)$ contains $P(\omega_1)$ as a 
subspace and each $U_n \cap P(\omega_1)$ is an open neighborhood of 
$\emptyset.$
\qed\enddemo

\proclaim{Corollary 4.10} Let $B$ be a complete Boolean algebra, and 
assume
that in the space $(B,\tau_s)$ there exists a countable family
$\{U_n\}_n$ of neighborhoods of $\boldkey 0$ such that $\bigcap_n
cl(U_n) =\{\boldkey 0\}.$ Then $B$ satisfies the countable chain 
condition
and $(B,\tau_s)$ is Fr\'echet.
\endproclaim

\demo{Proof}  $B$ satisfies ccc by Corollary 4.9. Also, $(B,\tau_s)$ 
is
clearly Hausdorff and so $B$ is weakly distributive by Lemma 4.2. 
Hence,
by Theorem 3.4, $B$ is a Fr\'echet space.
\qed\enddemo

\medpagebreak

We conclude the Section with some remarks:

A Fr\'echet space is Hausdorff if and only if
$$
\bigcap \{V\lor V : V\in \Cal N_0\} = \{0\}.
$$
Even more is true: If the space $(B,\tau_s)$ is Fr\'echet
and Hausdorff, then for every $k$,
$$
\bigcap\{V\lor V\lor ... \lor V (k \text{ times) }: V\in \Cal N_0\} = 
\{0\}.
$$
This is a consequence of the following:

\proclaim{Lemma 4.11} Let $B$ be a $\sigma-$complete Boolean algebra 
such that
$(B,\tau_s)$ is Fr\'echet. Then for every $U\in \Cal N^d_0$ there 
exists
a $V\in \Cal N_0^d$ such that $V\vee V\vee V\subseteq U\vee U.$
\endproclaim

In the next Section we use this consequence of Lemma 4.11:

\proclaim{Corollary 4.12} If $B$ is as in Lemma 4.11 and $U\in \Cal N^
d_0$
then there exists a $V\subseteq U$ in $\Cal N^d_0$ such that
$cl(V)\vee cl(V)\subseteq U\vee U.$
\endproclaim

(To see that this follows from Lemma 4.11, use $cl(V)\subseteq V\vee 
V.$)

\demo{Proof of Lemma 4.11}
Let us assume that for every $V\in \Cal N^d_0$ there exist $x,$ $y$ 
and $z$
in $V$ such that $x\vee y\vee z \notin U\vee U.$ Note that $U\vee U=
U\vartriangle U$ and is downward closed.

Let $V_0=U$; by induction we define neighborhoods $V_n$ and points
$x_n,y_n,z_n$ as follows: For each $n$ let
$x_n,y_n,z_n\in V_n$ be such that $x_n\vee y_n\vee z_n
\notin U\vee U.$ Then let $V_{n+1}\subseteq V_n$ be in 
$\Cal N^d_0$ such that the sets $x_n\vee V_{n+1},y_n\vee V_{n+1}$ and
$z_n\vee V_{n+1}$ are all included in $V_n;$ such a neighborhood 
exists
by the one sided continuity of $\vee.$

Let $X=\bigcap_n cl(V_n)$ and $\overline{x} = \overline{\lim_n}\,x_n,
\overline{y} = \overline{\lim_n}\,y_n, \overline{z} = 
\overline{\lim_n}\,z_n.$
The set $X$ is topologically closed and downward closed, and 
$X\subseteq cl(U)\subseteq U\vee U.$

We claim that $\overline{x},\overline{y},\overline{z}\in X.$ Thus let 
us
prove that for each $n,$ $\overline{x}\in cl(V_n).$ We have $x_n\in 
V_n,$
and by induction on $k>0$ we see that $x_n\vee x_{n+1}\vee...\vee 
x_{n+k}
\in V_n.$ Thus $\bigvee_{i\ge n}x_i \in cl(V_n)$ and $\overline{x}\in
cl(V_n).$

Next we claim that $\overline{x}\vee X\subseteq X$ (and similarly for
$\overline{y},\overline{z}).$ Let $n$ be arbitrary and let us show 
that
$\overline{x}\vee X \subseteq cl(V_n).$ If $k$ is arbitrary, we have
$x_n\vee...\vee x_{n+k}\vee V_{n+k+1}\subseteq V_n,$ and by the one 
sided
continuity of $\vee$ it follows that 
$x_n\vee...\vee x_{n+k}\vee cl(V_{n+k+1})\subseteq cl(V_n).$
Hence $x_n\vee...\vee x_{n+k}\vee X\subseteq cl(V_n),$ and so
$\bigvee_{i\ge n} x_i \vee X\subseteq cl(V_n).$
As $\overline{x}\le \bigvee_{i\ge n} x_i$ and $cl(V_n)$ is downward 
closed,
we have $\overline{x}\vee X\subseteq cl(V_n).$

Now it follows that $\overline{x}\vee\overline{y}\vee\overline{z}$ is 
in
$X$ and hence in $U\vee U.$ But $\overline{x}\vee\overline{y}\vee
\overline{z}=\overline{\lim_n}\,(x_n\vee y_n\vee z_n)$. As the 
complement
of $U\vee U$ is upward closed, we have 
$\bigvee_{i\ge n}(x_i\vee y_i\vee z_i)
\notin U\vee U$ for each $n$, and because $U\vee U$ is topologically 
open,
we have $\overline{x}\vee\overline{y}\vee\overline{z}\notin U,$ a 
contradiction. 
\qed\enddemo

\subheading {5. Metrizability}

We will show that for complete ccc Boolean algebras, Hausdorffness of the 
sequential topology is a strong property: it implies metrizability, and 
equivalently, the existence of a strictly positive Maharam submeasure. We remark 
that the assumption of completeness is essential.

\proclaim{Theorem 5.1} If $B$ is a complete Boolean algebra, then the 
following
are equivalent:
\roster
\item"(i)" $B$ is ccc and $(B,\tau_s)$ is a Hausdorff space,
\item"(ii)" there exists a countable family $\{U_n\}_n$ of open 
neighborhoods
of $\boldkey 0$ such that \linebreak $\bigcap_n cl(U_n)=\{\boldkey 0\},$
\item"(iii)" the operation $\lor$ is continuous at 
$(\boldkey 0,\boldkey 0)$, i.e. for every
$V\in\Cal N_0$ there exists a $U\in\Cal N_0$ such that $U\lor 
U\subseteq V,$
\item"(iv)" $(B,\tau_s)$ is a regular space,
\item"(v)" $(B,\tau_s)$ is a metrizable space,
\item"(vi)" $B$ carries a strictly positive Maharam submeasure.
\endroster
\endproclaim

The equivalence of (v) and (vi) is proved in \cite{Ma}, and (v) 
implies (i). We shall prove in this Section that 
properties (i)-(iv) are equivalent and imply (vi). First we claim that 
each
of the four properties implies that $B$ satisfies ccc, and that the 
space
$(B,\tau_s)$ is Fr\'echet.

If $B$ is ccc and Hausdorff, then by Theorems 4.1 and 3.4 it is 
Fr\'echet.

Property (ii) implies Fr\'echet by Corollary 4.10, and property (iv)
implies (i) (and hence Fr\'echet) by Corollary 4.8.

To complete the claim, 5.2--5.5 below prove that (iii) implies 
Fr\'echet.
Let $B$ be a complete Boolean
algebra and assume that $\vee$ is continuous at $(\boldkey 0,\boldkey 
0).$

\proclaim{Lemma 5.2} $B$ satisfies the countable chain condition.
\endproclaim

If $B$ does not satisfy ccc then $(B,\tau_s)$ contains $P(\omega_1)$ 
as
a closed subspace. Thus the Lemma 
is a consequence of the following lemma  closely related to Corollary 
4.7:

\proclaim{Lemma 5.3} In $(P(\omega_1),\tau_s)$ the operation $\cup$ is 
not
continuous at $(\emptyset,\emptyset).$
\endproclaim

\demo{Proof}
Let $U$ be the set of all $x\subset \omega_1$ whose complement is 
uncountable.
$U$ is an open neighborhood of $\emptyset$. 
We will show that for every
$V\in \Cal N_0$ there exist $Y$ and $Z$ in $V$ such that $Y\cup 
Z\notin U.$
Thus let $V\in \Cal N_0.$

By Lemma 4.6 there exists an $X$ such that
$X\notin U$ while  $X\in cl(V).$ By Corollary 3.7 there exist
$Y$ and $Z$ in $V$ such that $X=Y\vartriangle Z.$ But $Y\vartriangle Z
\subseteq Y\cup Z$ and therefore $Y\cup Z\notin U.$ 
\qed\enddemo

\proclaim{Lemma 5.4} $B$ is weakly distributive.
\endproclaim

\demo{Proof}
Assume that $B$ is not weakly distributive. By Lemma 4.2 there exists
some $a\ne 0$ such that $a\in cl(V)$ for every $V\in\Cal N_0.$

Let $U = \{x\in B: x\not\ge a\};$ $U$ is a neighborhood of $\boldkey 
0.$
We claim that for every $V\in\Cal N_0,$ $V\vee V\not\subseteq U,$ 
contradicting the continuity of $\vee.$ Thus let $V\in\Cal N_0$ be 
arbitrary.

We have $a\in cl(V).$ By Corollary 3.7, $a\in V\vartriangle V$ and so
there exist $x$ and $y$ in $V$ such that $a=x\vartriangle y.$ When we 
let
$b=x\vee y$ then $b\ge a$ and therefore $b\notin U$. But $b\in V\vee 
V,$
completing the proof.
\qed\enddemo

\proclaim{Corollary 5.5} $(B,\tau_s)$ is Fr\'echet.
\endproclaim

\demo{Proof} By Theorem 3.4.
\qed\enddemo

For the rest of section 5 we assume that $B$ is a complete Boolean 
algebra
that satisfies the countable chain condition, and that the space 
$(B,\tau_s)$
is Fr\'echet. In particular, $\Cal N_0^d$ is a neighborhood base, so 
we shall
only consider those neighborhoods of $\boldkey 0$ that are downward 
closed.

To prove that (i)-(iv) are equivalent, we first observe that (iii) 
implies (iv):

\proclaim{Proposition 5.6} If $\vee$ is continuous at 
$(\boldkey 0,\boldkey 0)$ then $(B,\tau_s)$ is regular.
\endproclaim

\demo{Proof} Let $V\in\Cal N_0.$ By homogeneity, it suffices to find
an open $U$ such that $cl(U)\subseteq V.$
Since $\vee$ is continuous at $(0,0)$ and 
since $(B,\tau_s)$ is Fr\'echet, there exists by Corollary 3.7 
a $U\in\Cal N^d_0$ such that $cl(U)\subseteq U\vee U\subseteq V.$
\qed\enddemo

As (iv) implies (i), it remains to show that (i) imples (ii) and that
(ii) implies (iii). Lemma 5.7 proves the latter:

\proclaim{Lemma 5.7} Assume that $(B,\tau_s)$ satisfies (ii).Then the
operation $\vee$ is continuous at $(\boldkey 0,\boldkey 0).$
\endproclaim

\demo{Proof} As $(B,\tau_s)$ is Fr\'echet, the set
$\Cal N^d_0$ of all downward closed open neighborhoods of $\boldkey 0$ 
is
a neighborhood base. Thus let us assume that there exists a $U\in \Cal 
N_0^d$
such that for every  $V\in \Cal N_0$ there exist $x$ and $y$ in $V$ 
with
$x\vee y \notin U.$

Let $\{V_n\}_n$ in $\Cal N_0^d$ be such that $\bigcap_n cl(V_n)=
\{\boldkey 0\}.$ We construct a descending sequence of neighborhoods 
$U_n$ in $\Cal N_0^d$
as follows: Let $U_0=V_0\cap U.$
Given $U_n$ let $x_n, y_n\in U_n$ be such that $x_n \vee y_n \notin 
U.$ By
(the separate) continuity of $\vee$ there exists a set 
$U_{n+1}\in \Cal N_0^d$ such that $x_n \vee U_{n+1}\subset U_n$ 
and $y_n \vee U_{n+1} \subset U_n;$ moreover, we may assume that 
$U_{n+1}$ is included in $V_{n+1}.$

Let $\overline{x}=\overline{\lim}\,x_n$ and $\overline{y}=
\overline{\lim}\,y_n.$ 
First we claim that $\overline{x}=\overline{y}=\boldkey 0$
and therefore $\overline{x} \vee\overline{y}=\boldkey 0 \in U.$

We have $\overline{x}=\bigwedge_n z_n$
where $z_n=\bigvee_k x_{n+k}.$ It suffices to prove that for each
$n$, $\bigwedge_m z_m$ is in the closure of $U_n$, and for that it is 
enough
to show that for each $m\ge n,$ $z_m\in cl(U_n).$ 

Let $n$ be arbitrary and ket $m\ge n$. As for each $k$ we have 
$x_{m+k} \vee U_{m+k+1}
\subset U_{m+k},$ it follows (by induction on $k$) that $x_m \vee 
x_{m+1}
\vee \dots \vee x_{m+k} \in U_m\subset U_n.$ Hence $z_m \in cl(U_n).$

Now we get a contradiction by showing that $\overline{x} \vee 
\overline{y}
\notin U.$ We have $\overline{x} \vee \overline{y} = \overline{\lim}
(x_n \vee y_n) = \bigwedge_n z_n$ where $z_n = \bigvee_{k\ge n} (x_k 
\vee y_k).$
As $U$ is a downward closed open set and $x_k \vee y_k \notin U$ for 
each $k,$ 
we have $z_n \notin U$ for each $n$ and therefore $\bigwedge_n z_n 
\notin U.$
\qed\enddemo

We now prove that (i) implies (ii):

\proclaim{Lemma 5.8} Let $B$ be a complete ccc Boolean algebra such 
that
$(B,\tau_s)$ is a Hausdorff space. Then there exists a sequence
$\{U_n\}_n$ in $\Cal N_0$ such that $\bigcap_n cl(U_n) =\{\boldkey 0\}.$
\endproclaim

\demo{Proof}
For any given $b\in B^+$ we shall find a sequence $\{V_n\}_n$ in $\Cal 
N_0^d$
such that $c_b=b-\bigvee(\bigcap_n cl(V_n))\ne 0.$ 
Then the set of all such $c_b$ is algebraically dense and therefore
there exists a partition $\{c_k\}_k$ of $\boldkey 1$ and sequences
$\{V_n^k\}_n$  with $\bigvee(\bigcap_n cl(V^k_n)) \wedge c_k=\boldkey 
0.$
Now when we let $U_n=V^0_n \cap V_n^1 \cap \dots \cap V_n^n$ for each 
$n$,
we get a sequence with the desired properties.

Thus let $b\ne \boldkey 0.$ We shall construct the sequence $\{V_n\}_n.$
For every set $S\subseteq B$ let $S^{(n)}$
denote the $n-$fold joint $S\vee...\vee S$ ($n$ times). 

As the space is Hausdorff, there exists a $V_0\in \Cal N^d_0$ such 
that
$b\notin V_0\vee V_0.$ By Lemma 4.11 and Corollary 4.12 there exists
for each $n$ some $V_{n+1}\in\Cal N^d_0$ such that $cl(V_{n+1})\vee
cl(V_{n+1})\subseteq V_n\vee V_n,$ and $V_{n+1}^{(3)}\subseteq V_n^
{(2)}.$
Let $X=\bigcap_n cl(V_n)$ and $a=\bigvee X.$

In order to prove that $b-a\ne\boldkey 0,$ it suffices to show that 
$a\in V_0\vee V_0,$ because that set is downward closed and $b$ is 
outside it.
By ccc, $a=\lim_n\,a_n$ where for each $n,$ $a_n\in X^{(n)}.$ We claim 
that
for each $n$, $X^{(n)}\subseteq V_2\vee V_2.$ Then $a\in cl(V_2\vee 
V_2)
\subseteq V_2^{(4)}\subseteq V_1^{(3)}\subseteq V_0^{(2)}.$

The claim is proved as follows (we may assume that $n$ is even):
$$
X^{(n)}\subseteq (cl(V_{n+1}))^{(n)}\subseteq V_n^{(n)}
\subseteq...\subseteq V_2^{(2)}.\qed
$$
\enddemo

This completes the proof of equivalence of properties (i)--(iv). 
We make the following remark:

\proclaim{Corollary 5.9} Let $B$ be a complete Boolean algebra such 
that
$(B,\tau_s)$ is a regular space. Then the he Boolean operations 
$\wedge,$ $-$ and 
$\vartriangle$ are continuous, and 
$(B,\vartriangle,\boldkey 0,\tau_s)$ is a
topological group. Moreover, $(B,\tau_s)$ is a completely regular 
space.
\endproclaim

\demo{Proof} As $(B,\tau_s)$ is Fr\'echet, $\boldkey 0$ has a 
neighborhood
base $\Cal N_0^d$ of sets for which $U\vartriangle U = U\vee U.$
Because $\vee$ is continuous at $(0,0)$, $\vartriangle$ is also 
continuous
at $(0,0)$. From that it easily follows that $\vartriangle$ is 
continuous
(at every $(u,v)\in B\times B$) and that $(B,\vartriangle,0,\tau_s)$ 
is
a topological group. Consequently, $\vee$ and $\wedge$ are also 
continuous
everywhere. 

Finally, every regular topological group is completely regular, cf. 
\cite{HeRo}.
\qed\enddemo

We shall now prove (vi), assuming that $(B,\tau_s)$ is regular.

\proclaim{Lemma 5.10} (a) There exists a sequence $\{U_n\}_n$ of 
elements of 
$\Cal N_0^d$ such that for every $n,$ $cl(U_{n+1})\subset U_{n+1} \vee
U_{n+1}\subset U_n$ and such that $\bigcap_n U_n = \{\boldkey 0\}.$ 

(b) Moreover, $\{U_n\}_n$ is a neighborhood base of $\,\boldkey 0.$
\endproclaim

\demo{Proof} (a) By continuity of $\vee$ there is a sequence $\{U_n\}_n$ 
in $\Cal N_0^d$ such that $U_{n+1} \vee U_{n+1} \subset U_n$ for every 
$n$.
By (ii) we may assume that $\bigcap_n U_n = \{\boldkey 0\}.$ 

(b) We prove that the $U_n$ form a 
neighborhood base. Assume not. Then there exists a $V\in \Cal N_0$ 
such that
for every $n,$ $U_n \not\subseteq V.$ For each $n$ let $x_n$ be such 
that
$x_n\in U_n-V.$

It follows by induction on $k$ that for each $n$ and each $k,$
$x_{n+1}\vee x_{n+2}\vee \dots \vee x_{n+k} \in U_n.$ Thus $\bigvee_k
x_{n+k}\in cl(U_n)$ and it follows that $\overline{\lim}\,x_n \in U_m$
for each $m;$ hence $\overline{\lim}\,x_n = \boldkey 0$ and so 
$\lim\,x_n=\boldkey 0.$ This is a contradiction because $V$ is a 
neighborhood of $\boldkey 0.$
\qed\enddemo

We are now ready to prove (vi). Let $\{U_n\}_n$ be a neighborhood base
of $\boldkey 0$ as in Lemma 5.10, with $U_0=B$.
Let $\bold D$ be the set of all rational
numbers of the form $r=\sum_{i=1}^{i=k} 2^{-n_i}$ where $\{n_1,\dots,
n_k\}$ is a finite increasing sequence of positive integers. For each
$r\in\bold D$ as above, let $V_r=U_{n_1}\vee\dots\vee U_{n_k},$ and 
let
$V_1=U_0=B.$
For each $a\in B,$ we define
$$
\mu(a) =\inf \{r\in \bold D \cup \{1\}: a\in V_r\}.
$$

\proclaim{Lemma 5.11} The function $\mu$ is a strictly positive 
Maharam 
submeasure.
\endproclaim

\demo{Proof} We repeatedly use the following fact, that follows by
induction on $k:$ For every increasing sequence $\{n_1,\dots,n_k\}$ of
nonnegative integers, $U_{n_1+1}\vee\dots\vee U_{n_k+1}\subseteq 
U_{n_1}.$

First, if $a\leq b$ then $\mu(a)\leq\mu(b);$ this is because for all
$r,s\in \bold D,$ if $r\leq s$ then $V_r\subseteq V_s.$

Second, for all $a$ and $b,$ $\mu(a\vee b)\leq \mu(a)+\mu(b);$ this is 
because
for all $r$ and $s$ such that $r+s<1,$ $V_r\vee V_s\subseteq V_{r+s}.$

Third, the submeasure $\mu$ is strictly positive: if $a\ne \boldkey 0$ 
then
there exists a positive integer $n$ such that $a\notin U_n=V_{1/2^n}$, 
and so
$\mu(a)\ge 1/2^n.$

Next we show that $\mu$ is continuous: if $\{a_n\}_n$ is a descending 
sequence
converging in $B$ to $\boldkey 0$ then for every $k$ eventually all
$a_n$ are in $U_k$, hence for eventually all $n$, $\mu(a_n)\leq 1/2^
k,$
and so $\lim_n\mu(a_n)=0.$

Finally, the topology induced by the submeasure $\mu$ coincides
with $\tau_s:$ this is because for each $n>0,$ $U_n\subseteq
\{a\in B: \mu(a)\leq 1/2^n\}\subseteq \bigcap_{k>n}(U_n\vee U_k)=
cl(U_n)\subseteq U_{n-1}.$
\qed\enddemo

\subheading{6. Sequential cardinals}

We now turn our attention to the atomic Boolean algebra $P(\kappa)$ 
where
$\kappa$ is an infinite cardinal. We compare two topologies on 
$P(\kappa):$
the product topology $\tau_c$ (when $P(\kappa)$ is identified with
the product space $\{0,1\}^\kappa$) and the sequential topology 
$\tau_s.$

If $f$ is a real-valued function on $B$ we
say that $f$ is {\it sequentially continuous} if it is continuous in 
the sequential topology $\tau_s$ on $B.$
Equivalently, $f(a_n)$ converges to $f(a)$ whenever $a_n$ converges
algebraically to $a$.

As $\tau_s$ is stronger than $\tau_c$, every real-valued function on 
$P(\kappa)$ that is continuous in the product topology is sequentially
continuous. Following \cite{AnCh} we say that $\kappa$ is a
{\it sequential cardinal} if there exists a discontinuous real valued
function that is sequentially continuous.

A submeasure $\mu$ on $P(\kappa)$ is {\it nontrivial} if $\mu(\kappa)>
0$ and
$\mu(\{\alpha\})=0$ for every $\alpha\in\kappa.$ If $\mu$ is a Maharam 
submeasure on $P(\kappa)$ then it is a sequentially continuous 
function.
If $\mu$ is nontrivial then it is discontinuous in the product 
topology,
because it takes value 0 on the dense set $[\kappa]^{\aleph_0}.$ Thus 
if
$P(\kappa)$ carries a nontrivial Maharam submeasure then $\kappa$ is
a sequential cardinal. In particular, the least real-valued measurable
cardinal is sequential. Keisler and Tarski asked in \cite{KeTa} 
whether
the least sequential cardinal is real valued measurable.

It follows from Theorem 6.2 below that 
if the Control Measure Problem has a positive answer then so does
the Keisler-Tarski question.

We use the following theorem of G. Plebanek (\cite{Pl}, Theorem 6.1).
A $\sigma-$complete Boolean algebra
$B$ carries a {\it Mazur functional} if there exists a sequentially
continuous real valued function $f$ on $B$ such that $f(\boldkey 0) =0$ and
$f(b) > 0$ for all $b\ne \boldkey 0.$

\proclaim{Theorem 6.1} {\rm (Plebanek)} If $\kappa$ is a sequential cardinal
then there exists a $\sigma$-complete proper ideal $H$
on $P(\kappa)$ containing all singletons such that the algebra $P(\kappa)/H$
carries a Mazur functional.
\endproclaim

\proclaim{Theorem 6.2} An infinite cardinal is sequential if and only if
the algebra $P(\kappa)$ carries a nontrivial Maharam submeasure.
\endproclaim

\demo{Proof} Let $\kappa$ be a sequential cardinal. By Theorem 6.1
the  $\sigma$-complete algebra $B=P(\kappa)/H$ carries a Mazur
functional $f$. First we claim that 
$B$ satisfies the countable chain condition, and hence is a complete 
algebra. If not, there is an uncountable antichain, and it follows 
that there is some $\varepsilon>0$ and
there are infinitely many pairwise disjoint elements $a_n,$ $n=0,1,2,...$
such that $|f(a_n)|\ge\varepsilon$ for all $n.$ This contradicts the 
sequential continuity of $f$ as $\lim_n \,a_n =\boldkey 0.$

For each $n,$ let
$U_n$ be the set of all $a\in B$ such that $|f(a)|<1/n.$ The $U_n$ are
neighborhoods of $\boldkey 0$ and satisfy property (ii) in Theorem 5.1.

By Theorem 5.1 $B$ carries a strictly positive  Maharam submeasure.
This submeasure induces a strictly positive Maharam submeasure on $P(\kappa)$
that vanishes $H$ and therefore on singletons. 
Thus $P(\kappa)$ carries a nontrivial Maharam submeasure.
\qed\enddemo

\enddocument
\end